\documentclass[12pt]{article}
\title{Equimultiple locus of embedded algebroid surfaces and blowing--up in
characteristic zero}
\author{Piedra-S\'anchez, R.\footnote{Supported by FQM 304 and BFM 2000--1523.}
 \and
Tornero, J.M.\footnote{Supported by FQM 218 and BFM 2001--3207.}}
\date{}

\newcommand{\NN}{{\bf N}}

\newcommand{\PP}{{\bf P}}

\newcommand{\vs}{\vspace{.15cm}}

\newcommand{\df}{\noindent {\bf Definition.-- }}

\newcommand{\thr}{\noindent {\bf Theorem.-- }}
\newcommand{\lema}{\noindent {\bf Lemma.-- }}
\newcommand{\dem}{{\bf Proof.-- }}
\newcommand{\obs}{\noindent {\bf Remark.-- }}

\newcommand{\ntt}{\noindent {\bf Notation.-- }}

\newcommand{\ord}{\mbox{ord}}
\newcommand{\Sp}{\mbox{Spec}}

\newcommand{\ol}{\overline}

\newcommand{\cE}{{\cal E}}
\newcommand{\cS}{{\cal S}}

\begin{document}

\maketitle

\abstract{The smooth equimultiple locus of embedded algebroid
surfaces appears naturally in many resolution process, both
classical and modern. In this paper we explore how it changes by
blowing--up.}

\vs


\vs


\section{Introduction}

During all this paper, we will consider $K$ an algebraically closed 
field of characteristic $0$ and $\cS = \Sp (K[[X,Y,Z]]/(F))$ an
embedded algebroid surface,
which, with no loss of generality, is considered to be defined
by a Weierstrass equation
$$
F(Z) = Z^n + \sum_{k=0}^{n-1} a_k (X,Y) Z^k,
$$
where $n$ is the multiplicity of $\cS$, that is $\ord (a_k) \geq
n-k$ for all $k=0,...,n-1$. After the well--known Tchirnhausen
transformation $Z \longmapsto Z - 1/a_{n-1}$ we can get a
Weierstrass equation of the form
$$
F(Z) = Z^n + \sum_{k=0}^{n-2} a_k (X,Y) Z^k.
$$

From now on, by a Weierstrass equation we will mean an equation
like this. Observe that, if we note the initial form of $F$ by $\ol{F}$,
the affine variety defined by $\ol{F}$ (that is, the tangent
cone of $\cS$) is a plane if and only if it is the plane $Z=0$.

In this situation the equimultiple locus of $\cS$ is
$$
\cE (\cS) = \left\{ P \in \Sp(K[[X,Y,Z]]) \; | \;
F \in P^{(n)} \right\},
$$
which is never empty, as $M=(X,Y,Z)$ always lies in $\cE (\cS)$.
Note that $Z \in P$ for all $P \in \cE(\cS)$.

Geometrically speaking, the equimultiple locus represents points
at where the multiplicity is the same than in the origin; hence
they are the ``closest'' points to the origin in (coarse) terms
of singularity complexity. We will note by $\cE_0 (\cS)$ the subset
of smooth elements of $\cE (\cS)$.

Our aim is studying the set $\cE(\cS)$ and, specifically, how its
elements change by blowing up, in order to have a better
understanding of the evolution of $\cE_0 (\cS)$ through a
resolution process. In this environment, our main result explains
how can we deduce ``geometrically'' $\cE_0 \left(
\cS^{(1)}\right)$ from $\cE_0 (\cS)$, where $\cS^{(1)}$ is the
result of blowing--up $\cS$ with center in an element of $\cE_0
(\cS)$.

The interest of this relies in the fact that the equimultiple
locus contains important information for desingularization
purposes. For instance, if $\cS$ has normal crossing
singularities, blowing--up centers lying in $\cE_0 (\cS)$ of maximal
dimension resolves the singularity at the origin. This is the
famous Levi--Zariski theorem on the resolution of surface
singularities, stated by Levi (\cite{Levi}) and proved by Zariski
(\cite{Z1}).

Concerning the extension of the results on this paper to the
arbitrary characteristic case, we must say that, no matter which
the characteristic is, the equimultiple locus has some very
interesting properties, some of which we state below:

(a) It is hyperplanar, that is, there exists a regular parameter
which lies in every element of $\cE(\cS)$. This was proved by
Mulay (\cite{Mulay}) after some previous work of Abhyankar
(\cite{GP}) and Narasimhan (\cite{Nar1}).

(b) As proved by Abhyankar (\cite{AbhB}), the Levi--Zariski
theorem remains true in positive characteristic.

However, the techniques used cannot be applied to the general
case, as it is not straightforward (actually, it is not known)
that the same regular parameter can define the tangent cone and a
hyperplane containing the equimultiple locus. In particular,
Mulay's proof of the hyperplanarity is not constructive, so a
different approach may be needed to handle with characteristic
$p>0$ (see \cite{RP} for some instances).

As for the extension of our result to higher dimensions is
concerned, we cannot be very optimistic. First of all, there is
the additional difficulty that, in positive characteristic, the
equimultiple locus is not hyperplanar anymore, as shown by
Narasimhan (\cite{Nar2}). Secondly, the Levi--Zariski resolution
process as stated does not work, as proved by Spivakovski
(\cite{Sp}), which, in any case, makes results on this line less
interesting. However, a less coarse version of the Levi--Zariski
theorem for a more general type of varieties will surely mean a 
great achievement and it will need results of this sort towards a 
more thorough understanding of the evolution of the equimultiple 
locus by sucessive blowing--ups.

\section{Notation and technical results}
\label{s:notation}

For the sake of completeness, we recall here basic facts and
technical results related to quadratic and monoidal
transformations that will be of some help in the sequel.

For all what follows, let $\cS$ be an embedded algebroid surface
of multiplicity $n$,
$$
F = Z^n + \sum_{k=0}^{n-2} \left( \sum_{i,j} a_{ijk}X^iY^j \right)
Z^k = Z^n + \sum_{k=0}^{n-2} a_k(X,Y) Z^k
$$
a Weierstrass equation of $\cS$. We will note
$$
N(F) = \left\{ (i,j,k) \in \NN^3 \; | \; a_{ijk} \neq 0
\right\}.
$$

\df
The elements of $\cE (\cS)$ different from $M$ will be called
equimultiple curves. The elements of $\cE_0(\cS)$ other than $M$ will be called
permitted curves.

\obs
The notion of permitted curves coincides with the
one derived from normal flatness in the work of Hironaka. 

\obs In particular, note that we can assume $P \in \cE_0 (\cS)$ to be,
for instance $(Z,X)$, after a suitable change of variables in
$K[[X,Y]]$. Plainly, $(Z,X)$ is permitted if and only if $i+k \geq
n$ for all $(i,j,k)\in N(F)$.

Under these circumstances, the monoidal transform of $\cS$,
centered in $(X,Z)$, in the point corresponding to the direction
$(\alpha:0 :\gamma)$ (say $\alpha \neq 0$) of the exceptional
divisor is the surface $\cS^{(1)}$ defined by the equation
$$
F^{(1)} = \left( Z_1 + \frac{\gamma}{\alpha}
\right)^n + \sum_{(i,j,k) \in N(F)} a_{ijk}X_1^{i+k-n} Y_1^j
\left( Z_1 + \frac{\gamma}{\alpha} \right)^k.
$$

Observe that this
only makes sense (that is, gives a non--unit) whenever
$\ol{F} (\alpha,0,\gamma)=0$.
The homomorphism
\begin{eqnarray*}
\pi^P_{(\ol{\alpha}:0:\gamma)}: K[[X,Y,Z]] & \longrightarrow &
K[[X_1,Y_1,Z_1]] \\
X & \longmapsto & X_1 \\
Y & \longmapsto & Y_1 \\
Z & \longmapsto & \displaystyle X_1 \left( Z_1 + \frac{\gamma}{\alpha}
\right) \\
\end{eqnarray*}
will be called the homomorphism associated to the monoidal
transformation in $(\ol{\alpha}:0:\gamma)$ or, in short,
the equations of the monoidal transformation. The overline is
because one must privilege a non-zero coordinate, but all the
possibilities define associated equations. 

As for quadratic transforms (that is, blowing--ups with center $M$)
is concerned: the quadratic
transform of $\cS$ in the point corresponding to the direction
$(\alpha:\beta:\gamma)$ (say $\alpha \neq 0$) of the exceptional
divisor is the surface $\cS^{(1)}$ defined by the equation
$$
F^{(1)} = \left( Z_1 + \frac{\gamma}{\alpha} \right)^n + \sum_{(i,j,k) \in N(F)} 
a_{ijk}X_1^{i+j+k-n} \left( Y_1 + \frac{\beta}{\alpha} \right)^j \left( Z_1 + 
\frac{\gamma}{\alpha} \right)^k.
$$

Again this only makes sense whenever
$\ol{F} (\alpha, \beta, \gamma)=0$. Analogously, the homomorphism
\begin{eqnarray*}
\pi^M_{(\ol{\alpha}:\beta:\gamma)}: K[[X,Y,Z]] & \longrightarrow &
K[[X_1,Y_1,Z_1]] \\
X & \longmapsto & X_1 \\
Y & \longmapsto & \displaystyle X_1 \left( Y_1 + \frac{\beta}{\alpha}
\right) \\
Z & \longmapsto & \displaystyle X_1 \left( Z_1 + \frac{\gamma}{\alpha}
\right) \\
\end{eqnarray*}
will be called the homomorphism associated to the quadratic
transformation in $(\ol{\alpha}:\beta:\gamma)$ or the equations of
the quadratic transformation. 

\obs
In the previous situation, consider a change of
variables in $K[[X,Y,Z]]$ given by
$$
\left\{ \begin{array}{ccc}
\varphi (X) & = & a_1 X' + a_2 Y' + a_3 Z' + \varphi_1 (X',Y',Z') \\
\varphi (Y) & = & b_1 X' + b_2 Y' + b_3 Z' + \varphi_2 (X',Y',Z') \\
\varphi (Z) & = & c_1 X' + c_2 Y' + c_3 Z' + \varphi_3 (X',Y',Z') \\
\end{array} \right.,
$$
with $\ord \left( \varphi_i \right) \geq 2$.

Assume also that
$$
\left\{ \begin{array}{ccc}
\alpha & = & a_1 \alpha' + a_2 \beta' + a_3 \gamma' \\
\beta & = & b_1 \alpha' + b_2 \beta' + b_3 \gamma' \\
\gamma & = & c_1 \alpha' + c_2 \beta' + c_3 \gamma' \\
\end{array} \right.
$$
with (say) $\gamma'\neq 0$. Then there is a unique
change of variables $\psi: K[[X_1,Y_1,Z_1]] \longrightarrow
K[[X'_1,Y'_1,Z'_1]]$ such that
$$
\psi \pi^M_{(\ol{\alpha}:\beta:\gamma)} = \pi^M_{(\alpha':\beta':
\ol{\gamma'})} \varphi.
$$

Both this remark and its monoidal counterpart (which will not be
used in this paper) are easy, although rather long, so we skip the
proofs. The interested reader may consult \cite{SG} and \cite{JMT}
for the complete details.

\df
Let $Q \in \cE(\cS)$, with $Q=(Z,G(X,Y))$. Then for
$\underline{u} \in \PP^2(K)$, the ideal
$$
\varpi^M_{\underline{u}} (Q) = \left( Z_1,
\frac{\pi^M_{\underline{u}}(G(X,Y))}{X_1^{\ord(G)}} \right)
$$
is called the quadratic transform of $Q$ in the point
corresponding to $\underline{u}$.

Obviously, this definition makes sense only if the quadratic
transform in the direction $\underline{u}$ does. There is a
natural version of monoidal transform of $Q$ with center $P$,
for all $P \in \cE_0(\cS)$.

\ntt
We will note by $\nu$ the natural isomorphism
$$
\nu : K[[X,Y,Z]]  \longrightarrow  K \left[\left[ X_1,Y_1,Z_1
\right]\right]
$$
sending $X$ to $X_1$, $Y$ to $Y_1$ and $Z$ to $Z_1$.

\section{The theorem}
\label{s:thr}

As our result is inspired by the resolution process, we will
restrict ourselves to the case which is interesting for
desingularization issues: that where $\cS$ and $\cS^{(1)}$ have
the same multiplicity. In particular, this avoids some
possibilities.

\lema
If the tangent cone of $\cS$ is not a plane, the multiplicity of any
monoidal transform of $\cS$ is strictly less than $n$.
 
\dem
If the tangent cone is not a plane, mind there is only one
possible element in $\cE_0 (\cS)$. After a suitable change of
variables on $K[[X,Y]]$, let this curve be $(Z,X)$. Then $\ol{F}$
cannot depend on $Y$,
$$
\ol{F} = Z^n + \sum_{i+k=n} a_{i0k} X^i Z^k = \prod_{l=1}^n
(Z-\alpha_l X).
$$

That is, the directions of the exceptional divisor are
$(1:0:\alpha_l)$ for $l=1,...,n$; not all of them equal.
Then, the equation for one of these monoidal transforms
is
$$
F^{(1)} = (Z_1+\alpha_{l_0})^n + \sum a_{ijk} X_1^{i+k-n}
Y_1^j (Z_1+\alpha_{l_0})^k
$$
where we find the monomial
$$
\left( \prod_{\alpha_{l_0} \neq \alpha_l}
\left( \alpha_{l_0} - \alpha_l \right) \right)
Z^m,
$$
with $m = \sharp\{ l \; | \; \alpha_{l_0} = \alpha_l \}$. This
monomial cannot cancel in any case. So there is a monomial in
$F^{(1)}$ of order strictly smaller than $n$ and we are done.

\thr
Let $\cS$ be an algebroid surface and $\cS^{(1)}$ a quadratic or 
monoidal transform of $\cS$ with the same multiplicity.

\noindent
{\bf (a)} If $\cS^{(1)}$ is the monoidal transform of $\cS$ with
center $P \in \cE_0 (\cS)$, then, either $\cE_0
\left( \cS^{(1)} \right) = \nu \left( \cE_0 (\cS) \right)$ or $\cE_0
\left( \cS^{(1)} \right) = \nu \left( \cE_0 (\cS) \setminus \{ P \}
\right)$.

\noindent
{\bf (b)} Let $\cS^{(1)}$ be the quadratic transform of $\cS$ in
the point corresponding to $\underline{u}$. Then:

{\bf (b.1)} If the tangent cone is not a plane then $\cE_0
\left( \cS^{(1)} \right) = \varpi^M_{\underline{u}} \left( \cE_0 (\cS) \right)$.

{\bf (b.2)} If the tangent cone is a plane, then in $\cE_0
\left( \cS^{(1)} \right)$ we can find three types of curves:

(i)  The exceptional divisor of the transform.

(ii) Primes $\varpi^M_{\underline{u}} (Q)$, with $Q \in \cE (\cS) \setminus
\cE_0(\cS)$, which are tangent to the exceptional divisor.

(iii) Primes $\varpi^M_{\underline{u}} (Q)$, with $Q \in \cE_0 (\cS)$, 
where both $\nu(Q)$ and $\varpi^M_{\underline{u}} (Q)$ are transversal 
to the exceptional divisor.

Moreover, if there is any prime of type (ii), it also appears the
type (i) prime.

\dem
We will do the proof case by case, although some arguments are
common to various instances. In what follows let $F$ be, as usual, a 
Weierstrass equation of $\cS$.

\vs

\noindent \fbox{{\bf Case (a)}}

\vs

From the previous lemma, we can assume that the tangent cone is
the plane $Z=0$. The basic tool for this situation is the
following:

\obs A monoidal transformation which does not imply a descent of
the multiplicity  cannot create new permitted curves.

This can  be proved easily as follows: 
under the hypothesis of case (a), let $(Z,G)$ be a permitted
curve. Then $\cS^{(1)}$, the monoidal transform of $\cS$ with
center at $(Z,G)$, is given by
$$
F^{(1)} = Z_1^n + \sum_{k=0}^{n-2} \frac{a_k(X_1,Y_1)}
{G(X_1,Y_1)^{n-k}} Z_1^k.
$$
 
This consists simply on taking $(Z,G)$ to $(Z,X)$ by a change
of variables (say $\varphi$) on $K[[X,Y]]$, applying the transform
(the only point in the exceptional divisor in this case is the
point corresponding to the direction $(1:0:0)$) and 
taking $\varphi^{-1}$ on $K[[X_1,Y_1]]$. The result
follows directly.

This remark clearly implies case (a) of the theorem.

\vs

\noindent \fbox{{\bf Case (b.2)}}

\vs

Some arguments in this case will be used for the other, so we will
begin for it. Let us start for the direction $(1:0:0)$ (the
direction $(0:1:0)$ is obviously symmetric). If
$$
F = Z^n + \sum_{(i,j,k) \in N(F)}
a_{ijk} X^iY^jZ^k,
$$
then
$$
F^{(1)} = Z_1^n + \sum_{(i,j,k) \in N(F)} a_{ijk}
X_1^{i+j+k-n}Y_1^jZ_1^k.
$$

Note that, as $F$ is a Weierstrass equation and the multiplicity
does not change, $F^{(1)}$ is a Weierstrass equation for
$\cS^{(1)}$, hence all the elements of $\cE \left( \cS^{(1)}
\right)$ are contained in $Z_1=0$, therefore all permitted curves
in $\cE \left( S^{(1)} \right)$ can be assumed to be of the form
$P = (Z_1, \gamma X_1 + \delta Y_1 + G(X_1,Y_1))$, with $\ord(G)
\geq 2$.

Let us prove now that, if a permitted curve transversal to the
exceptional divisor appears in $\cE_0 \left( \cS^{(1)} \right)$,
it comes from a permitted curve in $\cE_0 (\cS)$ which was also
transversal to the exceptional divisor (up to the action of $\nu$).

Suppose we have such a curve (that is, a prime as above with
$\delta \neq 0$). Then, applying the Weierstrass Preparation
Theorem, we may write $P$ as $(Z_1, Y_1+H(X_1))$. We have the
diagram

\unitlength=.9cm
\begin{picture}(11.5,5)(.5,-1)

\put(2,3){\makebox(0,0){$K[[X,Y,Z]]$}}
\put(10,3){\makebox(0,0){$K[[X',Y',Z']]$}}
\put(2,0){\makebox(0,0){$K[[X_1,Y_1,Z_1]]$}}
\put(10,0){\makebox(0,0){$K[[X'_1,Y'_1,Z'_1]]$}}

\put(3.5,3){\vector(1,0){5}}
\put(3.5,0){\vector(1,0){5}}
\put(2.5,2.5){\vector(0,-1){2}}
\put(9.5,2.5){\vector(0,-1){2}}

\put(6,3.25){\makebox(0,0){$\scriptstyle \varphi$}}
\put(6,0.25){\makebox(0,0){$\scriptstyle \psi$}}
\put(2.25,1.5){\makebox(0,0)[r]{$\scriptstyle \pi^M_{(1:0:0)}$}}
\put(9.75,1.5){\makebox(0,0)[l]{$\scriptstyle \pi^M_{(1:0:0)}$}}

\end{picture}

\noindent
with changes of variables
$$
\left\{ \begin{array}{ccl}
\varphi (X) & = & X' \\
\varphi (Y) & = & Y' - X'H(X') \\
\varphi (Z) & = & Z'
\end {array} \right. , \quad
\left\{ \begin{array}{ccl}
\psi (X_1) & = & X'_1 \\
\psi (Y_1) & = & Y'_1 - H(X'_1) \\
\psi (Z_1) & = & Z'_1
\end {array} \right.
$$

So, looking at the right vertical arrow, we have found a quadratic
transform on the direction $(1:0:0)$ which gives rise to the
permitted curve $(Z'_1,Y'_1)$. This clearly implies that $(Z',Y')$
was permitted in $\cS$. This proves the assertion.

Another way of seeing this is saying that, if there were no
permitted curves which were transversal to the exceptional
divisor, all permitted curves after the blowing--up must be
tangent to it.

Let us prove now that, if $(Z_1,X_1+Y_1^s v(Y_1))$ with $s>1$,
appears, so does $(Z_1,X_1)$. Write
$$
F^{(1)} = Z_1^n + \sum_{k=0}^{n-2} a^{(1)}_k (X_1,Y_1) Z_1^k,
$$
where it must hold $a^{(1)}_k = (X_1+Y_1^s v(Y_1))^{n-k}
b_k^{(1)}(X_1,Y_1)$.

Fix then $k_0 \in \{0,...,n-2\}$ and choose from all monomials in
$b_{k_0}^{(1)}$ the minimal one for the lexicographic ordering,
say $X_1^{i_0}Y_1^{j_0}$. Then all monomials appearing in
$a_{k_0}^{(1)}$ have exponent in $X_1$ greater or equal than $i_0$
and, besides, the monomial $X_1^{i_0} Y_1^{j_0+s(n-k_0)}$ actually
appears in $a_{k_0}^{(1)}$, as it cannot be cancelled.

Now it is plain that $(i,j,k) \in N \left( F^{(1)} \right)$
if and only if $( i-j-k+n, j, k) \in N(F)$, so
$$
i_0 \geq s\left(n-k_0\right)+k_0-n  \geq n-k_0.
$$

Hence $(Z_1,X_1) \in \cE \left( \cS^{(1)} \right)$.

Let us prove then the existence, in this case, of the curve $Q$ on
$\cE(\cS)$ announced in the theorem. As previously, we will
consider $\alpha=0$. We will prove that there exists a power series
$H(X,Y)$ verifying:
\begin{enumerate}
\item[(1)] $\ord (H) = \ord (G) = \lambda>1$.
\item[(2)] $H$ is regular on $Y$ of order $\lambda$.
\item[(3)] There is a unit $u(X_1,Y_1)$ such that
$$
\frac{1}{X_1^\lambda}H(X_1,X_1Y_1)=u(X_1,Y_1)(X_1+G(Y_1)).
$$
\end{enumerate}

This implies (quite straightforwardly) that $F \in Q^{(n)}$ and
$(Z,X+G(Y)) =\pi_{(1:0:0)}^M (Q)$, with $Q =(Z,H(X,Y))$. The second part 
is trivial and, for the first part it is enough proving
$$
\left( X_1+G(Y_1) \right)^{n-k} | a_k^{(1)} \left( X_1,Y_1 \right) \; 
\Longrightarrow \; H(X,Y)^{n-k} | a_k(X,Y),
$$
for $k=0,...,n-2$. Assume it is not so; then by the Weierstrass
Preparation Theorem and being $H$ regular with respect to $Y$,
we can write
$$
a_k (X,Y) = q(X,Y) H(X,Y)^{n-k} + \sum_{j=0}^{(n-k)\lambda - 1} 
\sigma_j (X) Y^j.
$$

Now we apply $\pi^M_{(1:0:0)}$ and we obtain
$$
\begin{array}{l}
X_1^a a^{(1)}_k \left( X_1,Y_1 \right)^{n-k} = \\
\quad \quad =
X_1^b q' \left( X_1,Y_1 \right) \left( X_1+G(Y_1) \right)^{n-k} + 
\sum_{j=0}^{(n-k)\lambda - 1} \left( \sigma_j (X_1) X_1^j \right) Y_1^j.
\end{array}
$$

Now, as $\left( X_1 + G ( Y_1) \right)^{n-k}$ divides $a_k^{(1)} \left( X_1,Y_1 \right)$
it also must divide $X_1^a a^{(1)}_k \left( X_1, Y_1 \right)$, hence the
uniqueness of quotient and remainder in the Weierstrass Preparation Theorem
imply 
$$
\sigma_j \left( X_1 \right) X_1^j = 0, \mbox{ for all } j=0,...,(n-k)\lambda - 1,
$$
and subsequently $H(X,Y)^{n-k} | a_k(X,Y)$. 

So let us prove the existence of $H$ and $u$. Write up $X_1+G(Y_1)$ as
$$
X_1 + G(Y_1) = X_1 + \sum_{i \geq \lambda} \alpha_i Y_1^i,
$$
and the power series we look for as
$$
H(X_1,Y_1) = \sum_{i+j=k \geq \lambda} \beta_{ij}X_1^iY_1^j, \; \;
u(X_1,Y_1) = \sum_{i+j=k \geq 0} \gamma_{ij} X_1^iY_1^j.
$$

It must hold
$$
\sum_{i+j=k\geq \lambda} \beta_{ij}X_1^{k-\lambda} Y_1^j =
\left(\sum_{i+j=k} \gamma_{ij} X_1^iY_1^j\right)
\left(X_1 + \sum_{i \geq \lambda} \alpha_i Y_1^i\right),
$$
which, for order 0, amounts to
$$
\beta_{\lambda, 0} = \gamma_{0,0}0 = 0.
$$

On the other hand, for order 1 we have
$$
\beta_{\lambda+1,0}X_1 + \beta_{\lambda-1,1}Y_1 =
\gamma_{0,0}X_1;
$$
that is, $\beta_{\lambda-1,1}=0$ and $\beta_{\lambda+1,0}=
\gamma_{0,0}$, whose value can be taken to be 1.

As for order 2,
$$
\beta_{\lambda+2,0}X_1^2 + \beta_{\lambda,1}X_1Y_1+
\beta_{\lambda-2,2}Y_1^2 = \gamma_{0,0}\alpha_2Y_1^2 + \gamma_{1,0}
X_1^2 + \gamma_{0,1}X_1Y_1,
$$
which forces $\beta_{\lambda-2,2}=\alpha_2$ and let us freedom
for fixing $\beta_{\lambda+2,0}=\gamma_{1,0}$, and
$\beta_{\lambda,1} = \gamma_{0,1}$.

Observe then the following facts:
\begin{itemize}
\item Each $\beta_{ij}$ appears only for order $i+2j-\lambda$.
\item In order $k$, all $\gamma_{ij}$ with $i+j<k$ appear,
but they never have relations among them, only those of the type
$$
\beta_{ab} = \sum \gamma_{cd}\alpha_e.
$$
\end{itemize}

Therefore it is clear that we can choose arbitrarily the
$\gamma_{ij}$, and this choosing determines the $\beta_{ij}$. Therefore
both $H$ and $u$ exist.

It only remains proving that $H$ can be chosen such that $H(0,Y)$ has
order $\lambda$. But this is direct from the formula for
order $\lambda$;
$$
\beta_{0,\lambda} Y_1^\lambda =  \gamma_{0,0} \alpha_\lambda Y_1^\lambda
\neq 0,
$$
so $\beta_{0,\lambda}\neq 0$.

For the results at points $(1:\alpha:0)$ it suffices considering
the (commutative) diagram

\unitlength=.9cm
\begin{picture}(11.5,5)(.5,-1)

\put(2,3){\makebox(0,0){$K[[X,Y,Z]]$}}
\put(10,3){\makebox(0,0){$K[[X',Y',Z']]$}}
\put(2,0){\makebox(0,0){$K[[X_1,Y_1,Z_1]]$}}
\put(10,0){\makebox(0,0){$K[[X'_1,Y'_1,Z'_1]]$}}

\put(3.5,3){\vector(1,0){5}}
\put(3.5,0){\vector(1,0){5}}
\put(2.5,2.5){\vector(0,-1){2}}
\put(9.5,2.5){\vector(0,-1){2}}

\put(6,3.25){\makebox(0,0){$\scriptstyle \varphi$}}
\put(6,0.25){\makebox(0,0){$\scriptstyle \psi$}}
\put(2.25,1.5){\makebox(0,0)[r]{$\scriptstyle \pi^M_{(1:\alpha:0)}$}}
\put(9.75,1.5){\makebox(0,0)[l]{$\scriptstyle \pi^M_{(1:0:0)}$}}

\end{picture}

\noindent
with $\varphi$ given by
$$
\left\{ \begin{array}{ccl}
\varphi (X) & = & X' \\
\varphi (Y) & = & Y' - \alpha X' \\
\varphi (Z) & = & Z'
\end {array} \right.
$$

\obs
This case, in geometrical terms, may be expressed as follows:
\begin{itemize}
\item Permitted curves transversal to the exceptional divisor cannot
be created nor erased.
\item Permitted curves tangent to the exceptional divisor are erased.
\item Permitted curves tangent to the exceptional divisor can be
created from desingularization of equimultiple (singular) curves. In 
this case, one of them must be the exceptional divisor itself.
\end{itemize}

\vs

\noindent \fbox{{\bf Case (b.1)}}

\vs

\obs
In the conditions of (b.1), let $P = (\alpha:\beta:\gamma)$ a
point in the tangent cone with multiplicity $r$. Then, the
quadratic transform of $\cS$ on $(\alpha:\beta:\gamma)$ has, at
most, multiplicity $r$. This is straightforward, using, for
instance, the Taylor expansion of $F$.

So we only need to be concerned about points of multiplicity $n$
on the tangent cone. Changing the variables if needed we may
consider that the point is $(0:1:0)$ and, subsequently, $\ol{F}$
does not depend on $Y$.

We will first prove that the quadratic transform cannot have
permitted curves. Note that, in (b.2), we have showed that, if
a new permitted curve appears, so does the exceptional divisor
(and we did not use the fact that $\ol{F} = Z^n$ for proving
this). But $(Z,Y)$ cannot be a permitted curve, $\ol{F^{(1)}}$
having monomials in $K[X,Z]$ other than $Z^n$.

Now we explain why the quadratic transform does not erase
permitted curves either. In fact if there is a permitted curve
(only one is possible), we may take it to be $(Z,X)$, after the
customary change of variables. Then it is plain that, after a
quadratic transform on the direction $(0:1:0)$, $(Z_1,X_1)$
remains permitted, simply looking at the characterization given in
the previous section.

This finishes the proof of the theorem.

\vs\vs

\noindent
\begin{tabular}{ll}
R. Piedra: & 
Departamento de \'Algebra, Universidad de Sevilla. \\
& Apdo 1160, 41080 Sevilla (Spain). \\
& E--mail: piedra@algebra.us.es

\\ \\

\noindent
J.M. Tornero: 
& Departamento de \'Algebra, Universidad de Sevilla. \\
& Apdo 1160, 41080 Sevilla (Spain). \\
& E--mail: tornero@algebra.us.es
\end{tabular}
\end{document}